\documentclass[12pt,a4paper,oneside]{amsart}

\usepackage{amsfonts, amsmath, amssymb, amsthm, amscd}
\usepackage{hyperref}
\hypersetup{
  colorlinks   = true, 
  urlcolor     = blue, 
  linkcolor    = blue, 
  citecolor   = red 
}
\usepackage{anysize}
\marginsize{2cm}{2cm}{1.5cm}{1.5cm}
\usepackage[pdftex]{graphicx}

\newtheorem{theorem}{Theorem}[section]
\newtheorem{lemma}[theorem]{Lemma}

\newtheorem{conjecture}[theorem]{Conjecture}
\newtheorem{corollary}[theorem]{Corollary}

\theoremstyle{definition}
\newtheorem{definition}[theorem]{Definition}

\theoremstyle{remark}
\newtheorem{remark}[theorem]{Remark}
\newtheorem{question}[theorem]{Question}

\newcounter{fig}

\numberwithin{equation}{section}

\newcommand{\pt}{\mathrm{pt}}
\newcommand{\cyc}{\mathit{cl}}
\newcommand{\CAT}{\mathbf{CAT}}

\DeclareMathOperator{\hind}{ind}

\DeclareMathOperator{\vol}{vol}
\DeclareMathOperator{\diam}{diam}
\DeclareMathOperator{\dist}{dist}

\DeclareMathOperator{\st}{st}

\renewcommand{\epsilon}{\varepsilon}

\renewcommand{\phi}{\varphi}

\title{Borsuk--Ulam type theorems for metric spaces}

\author{Arseniy~Akopyan{$^\spadesuit$}}
\email{akopjan@gmail.com}

\author{Roman~Karasev{$^\clubsuit$}}
\email{r\_n\_karasev@mail.ru}
\urladdr{http://www.rkarasev.ru/en/}

\author{Alexey~Volovikov{$^\diamondsuit$}}
\email{a\_volov@list.ru}

\thanks{{$^\spadesuit$} Supported by the European Research Council (ERC) under the European Union's Horizon 2020 research and innovation programme (grant agreement  No 716117)}

\address{{$^\spadesuit$} Institute of Science and Technology Austria (IST Austria), Am Campus 1, 3400 Klosterneuburg, Austria}
\address{{$^\clubsuit$} Institute for Information Transmission Problems RAS, Bolshoy Karetny per. 19, Moscow, Russia 127994 and Moscow Institute of Physics and Technology, Institutskiy per. 9, Dolgoprudny, Russia 141700 }

\address{{$^\diamondsuit$} Department of Higher Mathematics, Moscow State Institute of Radio-Engineering, Electronics and Automation (Technical University), Pr. Vernadskogo 78, Moscow 117454, Russia}


\begin{document}

\begin{abstract}
In this paper we study the problems of the following kind: For a pair of topological spaces $X$ and $Y$ find sufficient conditions that under every continuous map $f : X\to Y$ a pair of sufficiently distant points is mapped to a single point.
\end{abstract}

\maketitle

\section{Introduction}

In this paper we are going to give new proofs, using the recent ideas of M.~Gromov, to the classical Borsuk--Ulam and Hopf theorems and their generalizations, and study some their consequences, as well as separate results, in the spirit of Urysohn width and Gromov waist.

Recall the famous Borsuk--Ulam theorem~\cite{bor1933}:

\begin{theorem}[K.~Borsuk, S.~Ulam, 1933]
\label{bu}
Under any continuous map $f : S^n\to \mathbb R^n$ some two opposite points are mapped to a single point.
\end{theorem}

A deep generalization of this result is the Hopf theorem~\cite{hopf1944}:

\begin{theorem}[H.~Hopf, 1944]
\label{theorem:hopf}
Let $X$ be a closed Riemannian manifold of dimension $n$ and $f : X\to \mathbb R^n$ be a continuous map. For any prescribed $\delta > 0$, there exists a pair $x,y\in X$ such that $f(x) = f(y)$ and the points $x$ and $y$ are connected by a geodesic of length $\delta$.
\end{theorem}



Our presentation is greatly inspired by the results of~\cite{gromov2010}, where the estimates for the size of the preimage of a point were proved using the technique of ``contracting in the space of (co)cycles''. One of the questions addressed in this paper is how this technique can be applied to the Borsuk--Ulam and Hopf theorems. Such an application turns out to be possible and these old theorems are generalized (see Theorems~\ref{theorem:bu-paths} and Theorem~\ref{theorem:hopf-paths}).

\textbf{Acknowledgment.}
The authors thank Sergey~Avvakumov, Alexey~Balitskiy, Vladimir Dol'nikov, Misha Gromov, Bernhard Hanke, and Evgeniy~Shchepin for useful discussions and remarks.

\section{A Borsuk--Ulam type theorem for metric spaces}
\label{cycl-sec}

We are going to utilize the ideas of M.~Gromov~\cite{gromov2010} to give a coincidence theorem. Let us make a few definitions.

\begin{definition}
Let $X$ be a topological space. Denote $PX$ the \emph{space of paths}, i.e. the continuous maps $c :[0,1]\to X$. This space has a natural $\mathbb Z_2$-action by the change of parameter $t\mapsto 1-t$, and a natural $\mathbb Z_2$-equivariant map
$$
\pi : PX \to X\times X,\quad c \mapsto (c(0), c(1)).
$$
\end{definition}

\begin{definition}
Call a $\mathbb Z_2$-equivariant section $s$ of the bundle $\pi :PX \to X\times X$ over an open neighborhood $\mathcal D(s)$ of the diagonal $\Delta(X)\subset X\times X$ a \emph{short path map}, iff $s(x, x)$ is a constant path for any $x\in X$.
\end{definition}

Such short path maps may be given by assigning a shortest path to a pair of points in a metric space. If $X$ is a compact Riemannian manifold then such short path maps do exist.

Now we are ready to state:

\begin{theorem}
\label{theorem:bu-paths}
Suppose $X$ is a closed manifold of dimension $n$, $Y$ is another manifold of dimension $n$, and $f : X\to Y$ is a continuous map of even degree. Then for any short path map $s : X\times X\to PX$ there exists a pair $(x,y)\not\in \mathcal D(s)$ such that $f(x) = f(y)$.
\end{theorem}

The classical Borsuk--Ulam theorem~\cite{bor1933} follows from this theorem by considering $X=S^n$ and $s$ to be the shortest path map in the standard metric. Theorem~\ref{theorem:hopf} (of Hopf) does not follow from this theorem because here we may only obtain an inequality on $\dist(x, y)$. An advantage of Theorem~\ref{theorem:bu-paths} is that the codomain $Y$ may be arbitrary.

In a similar way as the Borsuk--Ulam theorem produces the ham sandwich theorem~\cite{stonetukey1942,steinhaus1945}, it is possible to produce a ham-sandwich type result from the above theorem:

\begin{theorem}
\label{hs-balls}
Let $M$ be a closed manifold of dimension $n$ with a path-metric and convexity radius $\kappa$. Assume we have $n$ Borel measures $\mu_1,\ldots,\mu_n$ in $M$ that are zero on every metric sphere in $M$ and a number $0\le r \le \kappa$. Then there exist two points $x,y\in M$ such that their neighborhoods $B_r(x)$ and $B_r(y)$ do not overlap and 
\[
\mu_i B_r(x) = \mu_i B_r(y)
\]
for every $i=1,\ldots, n$. 
\end{theorem}

The original ham-sandwich follows from this theorem if we put $M = S^n$ and $r = \pi/2$.

\begin{proof}
The assumptions on the measure guarantee that the values
\[
\mu_1 B_r(x) , \ldots, \mu_n B_r(x)
\]
produce a continuous map $f : M \to \mathbb R^n$. Applying Theorem~\ref{theorem:bu-paths} to $f$ and the shortest paths in $M$, we obtain two points $x,y\in M$ with $\mu_i B_r(x) = \mu_i B_r(y)$ for every $i=1,\ldots, n$. If the balls $B_r(x)$ and $B_r(y)$ were overlapping (had a common interior point) then the shortest path between $x$ and $y$ would be unique (see Remark~\ref{conv-rad-remark} below), a contradiction.
\end{proof}

\begin{definition}
Suppose $X$ is a compact Riemannian manifold. Let $\kappa(X)$ be the maximum number such that for any $0<\delta<\kappa(X)$ any ball in $X$ of radius $\delta$ is geodesically strictly convex. Call $\kappa(X)$ the \emph{convexity radius}. 
\end{definition}

\begin{remark}
\label{conv-rad-remark}
Obviously $\rho(X) \ge 2 \kappa(X)$, because touching strictly convex balls can intersect at one point only. It is also known that $\kappa(X)>0$ for compact Riemannian manifolds.
\end{remark}

\section{Space of cycles and the proof of Theorem~\ref{theorem:bu-paths}}

We start from explaining the main ideas underlying what Gromov calls ``contraction in the space of cycles''~\cite{gromov2010} in a particular case. Denote $\cyc_0(X; \mathbb F_2)$ the space of $0$-cycles mod $2$ in $X$, that is the space of formal finite combinations $\sum_{x\in X} a_x x$ with $a_x\in\mathbb F_2$ and $\sum_{x\in X} a_x = 0$ with an appropriate topology. 

A more tangible description of $\cyc_0(X;\mathbb F_2)$ is the union over $k\ge 0$ of spaces of unordered $2k$-tuples $B(X, 2k) \subset X^{\times 2k} / \Sigma_{2k}$. Informally, the topology in $\cyc_0(X;\mathbb F_2) = \bigcup_{k\ge 0} B(X, 2k)$ is such that when two points of a set $c\in B(X, 2k)$ tend to a single point then they ``annihilate'' giving a configuration in $B(X, 2k-2)$ in an obvious way, and conversely a pair of points may be ``created'' from a single point giving a configuration in $B(X, 2k+2)$.

In the case when $X$ is an $n$-dimensional manifold we define the \emph{canonical class} $\xi$ in $H^n(\cyc_0(X;\mathbb F_2); \mathbb F_2)$ as follows. Any $n$-dimensional homology class of $\cyc_0(X; \mathbb F_2)$ can be represented by a chain $c$, which is given by a map of an $n$-dimensional mod $2$ pseudomanifold $K$ to $\cyc_0(X;\mathbb F_2)$. Considering any element of $\cyc_0(X;\mathbb F_2)$ as a subset of $X$ we may consider $c$ as a set valued map from $K$ to $X$. Its graph $\Gamma_c$ is a subset of $K\times X$, which is again a mod $2$ pseudomanifold, and the projection $\Gamma_c\to K$ has degree zero. Hence the degree mod $2$ of the natural projection $\Gamma_c \to X$ is well defined. This degree will be the value $\xi(c)$ by definition. Another informal way to define $\xi$ would be to count how many times a generic point $x_0\in X$ participates in the $0$-cycles from the chain $c$.

Now we return to the proof of the theorem. From the compactness considerations it is sufficient to prove the theorem for smooth generic maps $f$. In this case we may define the natural map 
$$
f^c : Y \to \cyc_0(X; \mathbb F_2),
$$
which maps any $y\in Y$ to the $\mathbb F_2$-cycle 
$$
f^c(y) = \sum_{x\in f^{-1}(y)} c_x x,
$$
where $c(x)$ is the mod $2$ multiplicity of the map $f$ at $x$. This map is well-defined because the degree of $f$ is even by the hypothesis. The image of $f^c$ represents an $n$-dimensional mod $2$ homology class in $\cyc_0(X;\mathbb F_2)$ and by the definition of the fundamental class $\xi$ it is obvious that $\xi(f^c(Y)) = 1$. Therefore the map $f^c$ is homotopically nontrivial.

But we are going to deform the map $f^c$ to the constant map by a homotopy $h_t$, using the short path map $s$. Put
$$
h_t(y) = \sum_{x_1\neq x_2\in f^{-1}(y),\ c_{x_1}, c_{x_2} = 1} s(x_1, x_2)(t/2). 
$$

We have to check whether this map is continuous in $y$ and $t$. If the preimage $f^{-1}(y)$ changes by ``annihilating'' a pair points or ``creating'' a pair of points, the components of $h_t(y)$ are also ``annihilated'' or ``created'' pairwise (here we use the $\mathbb Z_2$-equivariance of the short path map and its behavior over the diagonal).

Let us explain the words ``created'' and ``annihilated''. For generic smooth $f : X\to Y$, let $S_1\subset Y$ be the set of special values of $f$, which has codimension at least $1$. Let the set $S_2\subset S_1$ correspond to the singularities of $f$ more complicated than folds. For generic smooth $f$, the set $S_2$ has codimension at least $2$ in $Y$ and its preimage $f^{-1}(S_2)$ has codimension at least $2$ in $X$. We may ignore $S_2$ in the reasoning with the fundamental class of $X$ or $Y$, because the homology is not affected by codimension $2$ changes. The space $Y\setminus S_2$ remains connected from dimension considerations whenever $Y$ was connected. Now we see that when the point $y$ travels in $Y\setminus S_2$, the graph $G_y$ may only change when $y$ crosses a fold singularity and some two vertices of the graph are ``created'' or ``annihilated''.

If the parameter $t$ approaches $0$ then $h_t$ approaches $f^c$, because in every $f^{-1}(y)$ we have an even number of points with odd multiplicities $c_x$, so in the expression of $h_t$ we approach every point $x\in f^{-1}(y)$ (such that $c_x$ is odd) odd number of times. If $t$ approaches $1$ the points $s(x_1, x_2)(t/2)$ and $s(x_2, x_1)(t/2)$ tend to ``annihilate'' (and do ``annihilate'' at $t=1$), and therefore $h_1$ maps the whole $Y$ to zero cycle. Thus the proof is complete.

\begin{remark}
\label{fund-filling}
In~\cite{kar2011} a simplified version of the reasoning in~\cite{gromov2010}, in the particular case of the problem of probability of covering by a simplex, was presented, which avoids an explicit use of the space of (co)cycles. In the above proof a similar trick is also possible in the following way. 

Assume that the map $f : X\to Y$ is generic in a certain sense. For example, when $X$ is triangulated and $Y$ is $\mathbb R^n$ then $f$ may be thought of as a generic PL map. Then for any $y\in Y$ consider the finite set $f^{-1}(y)$ and the complete graph ($1$-dimensional complex) $G_y$ on the vertices $f^{-1}(y)$. Denote the union of these complete graphs over $y\in Y$ by $G_f$. With some natural topology (starting with the topology of $X\times Y$) $G_f$ can be interpreted as an abstract chain (in PL case this can be made rigorous by endowing $G_f$ a CW structure). 

The boundary of $G_f$ module $2$ and modulo codimension $2$ is not generated by the only condimension $2$ singularity of $f$, the fold (in PL case a fold is the situation when the two top dimensional faces $\sigma$ and $\tau$ are mapped to the one side of the image of there common codimension $1$ face $\rho$), which roughly corresponds to what is called ``creation'' and ``annihilation'' above. The remaining part of the boundary of $G_f$ modulo $2$ is $\partial G_f = \bigcup_{x\in X} \deg_{G_{f(x)}} x$. Under the assumption that the degree of $f$ is even, we conclude that generically a vertex of $G_y$ has odd degree, and therefore $\partial G_f = \bigcup_{x\in X} x$ modulo $2$.

If any edge of $G_y\subset G_f$ can be realized in $X$ with continuous dependence on the endpoints (for example, using a short path map) then $G_f$ is continuously mapped to $X$. So $G_f$ becomes an $(n+1)$-dimensional chain in $C_{n+1}(X; \mathbb F_2)$ with boundary $\partial G_f = [X]$ modulo $2$. But the fundamental class of a closed manifold $X$ modulo $2$ cannot vanish, which is a contradiction.
\end{remark}

\section{Classical Hopf type results}
\label{section:hopf-classical}

First, for completeness, we remind the proof of the Hopf theorem (Theorem~\ref{theorem:hopf}). The proof is given in the original paper~\cite{hopf1944} in German and it makes sense to repeat it here in English.

\begin{proof}[Proof of Theorem~\ref{theorem:hopf}, translated from its German version in~\cite{hopf1944}]
Choose a point $p\in f(X)$ that has the maximal coordinate $x_1$ among the image $f(X)$. Let $o$ be any point in $f^{-1}(p)$.

Let $T$ be the tangent space of $X$ at $o$, $S$ be its unit sphere, and $\phi : T \to X$ be the exponential map. We define a continuous family of maps
\[
h_t : S\to \mathbb R^n,\quad h_t(v) = f\left(\phi\left(\frac{\delta(t+1)}{2} v\right)\right) - f\left(\phi\left(\frac{\delta(t-1)}{2} v\right)\right).
\]
Note that the points $x = \phi\left(\frac{\delta(t+1)}{2} v\right)$ and $y=\phi\left(\frac{\delta(t-1)}{2} v\right)$ are always connected by a geodesic of length $\delta$, which is the exponential image of a straight line segment through the origin in $T$. If we assume that $f(x)$ is never equal to $f(y)$, then $h_t(v)$ is never zero and we can define
\[
\tilde h_t(v) = \frac{h_t(v)}{|h_t(v)|}.
\]
Now we note that $\tilde h_0$ is an odd map between $(n-1)$-dimensional spheres and therefore has odd degree (this is a consequence of the Borsuk--Ulam theorem). By the choice of $p=f(o)=f(\phi(0))$ with maximal first coordinate we have that the first coordinate of 
\[
h_1(v) = f\left(\phi\left( \delta v\right)\right) - f\left(\phi\left(0\right)\right)
\]  
is not positive. So the map $\tilde h_1$ does not contain the vector $(1,0,\ldots, 0)$ in its image and therefore must have zero degree. Now we obtain a contradiction because $\tilde h_0$ is homotopic to $\tilde h_1$ and the degree of a map is a homotopy invariant. Hence $f(x)=f(y)$ for some pair $(x, y)$ connected by a geodesic of length $\delta$.
\end{proof}

A certain extension of this result is known. Here we provide the statement and a proof of a particular case of the result in~\cite{tomdiecksmith1979}, avoiding the use of localization techniques and invoking the Adams theorem~\cite{ad1960} instead in the proof that we provide here.

\begin{theorem}[T.~tom Dieck, L.~Smith, 1979]
\label{theorem:hopf-degree-zero}
Let $n$ be a positive integer not equal to $1$, $3$, or $7$, and let $f : S^n\to S^n$ be a continuous map of even degree. For any prescribed $\delta>0$ and any Riemannian metric on the sphere $S^n$, there exists a pair $x,y\in S^n$ such that $f(x) = f(y)$ and the points $x$ and $y$ are connected by a geodesic of length $\delta$.
\end{theorem}

\begin{proof}
Let $U$ be the unit tangent vector bundle of $S^n$, that is the set of pairs $(x,v)$ of unit vectors in $\mathbb R^n$ such that $x\cdot v = 0$. This space $U$ has an involution $(x,v)\mapsto (x,-v)$. Assuming that the problem has no solution for a map $f$ and a number $\delta$ we build an equivariant map $F : U\to S^n$, where the involution on $S^n$ is $y\mapsto -y$.

Take the geodesic $\gamma_{x,v}$ on the sphere such that $\gamma_{x,v}(0) = x$ and $\dot\gamma_{x,v}(0)=v$. Then the two points $\gamma_{x,v}(-\delta/2)$ and $\gamma_{x,v}(\delta/2)$ are connected by a geodesic of length $\delta$ and this allows to correctly define
\[
F(x,v) = \frac{f(\gamma_{x,v}(\delta/2)) - f(\gamma_{x,v}(-\delta/2))}{\left| f(\gamma_{x,v}(\delta/2)) - f(\gamma_{x,v}(-\delta/2)) \right|}.
\]
Changing $v$ to $-v$ interchanges the points $\gamma_{x,v}(-\delta/2)$ and $\gamma_{x,v}(\delta/2)$ and shows that $F(x,-v) = - F(x,v)$, that is $F$ is equivariant. 

For another variable $t\in [0,\delta/2]$, one considers the homotopy
\[
\frac{f(\gamma_{x,v}(\delta/2 - t)) - f(\gamma_{x,v}(-\delta/2-t))}{\left| f(\gamma_{x,v}(\delta/2-t)) - f(\gamma_{x,v}(-\delta/2-t)) \right|},
\]
which is well-defined because the two points $\gamma_{x,v}(-\delta/2-t)$ and $\gamma_{x,v}(\delta/2-t)$ are connected by a geodesic of length $\delta$ and are not mapped to a single point by $f$. This homotopy (non-equivariantly) connects $F$ to the map
\[
G(x,v) = \frac{f(\gamma_{x,v}(0)) - f(\gamma_{x,v}(\delta))}{\left| f(\gamma_{x,v}(0)) - f(\gamma_{x,v}(\delta)) \right|} = \frac{f(x) - f(\gamma_{x,v}(\delta))}{\left| f(x) - f(\gamma_{x,v}(\delta)) \right|},
\]

Let $D$ be the tangent disk bundle of $S^n$, whose boundary is $U$. The map $G: U\to S^n$ trivially extends to $D$ by 
\[
\frac{f(x) - |v| f(\gamma_{x,v/|v|}(\delta))}{\left| f(x) - |v| f(\gamma_{x,v/|v|}(\delta)) \right|}.
\]
The homotopy extension theorem then implies that the equivariant map $F : U\to S^n$ gets non-equivariantly extended to $D\to S^n$. Let $T$ be the unit sphere bundle of $TS^n\oplus \tau$ where $\tau$ is the trivial one-dimensional bundle. The space $T$ has the splitting in two copies of $D$ corresponding to the positive and negative direction of $\tau$. The extension of $F$ to the positive copy of $D$ can be equivariantly extended to its negative copy to yield a continuous equivariant map $\widetilde F : T \to S^n$. Since $\tau$ may be considered as the normal bundle of $S^n$ in $\mathbb R^{n+1}$, we actually have $TS^n\oplus \tau\cong \tau^{n+1}$ and $T\cong S^n\times S^n$. The equivariance of the map $\widetilde F : S^n\times S^n\to S^n$ is then understood so that $\widetilde F(-x,y)=-\widetilde F(x,y)$.

From the definition of $G$ and $\widetilde F$ it also follows that the composition of the diagonal inclusion $S^n\to S^n\times S^n$, $x\mapsto (x,x)$, and $\widetilde F$ is homotopic to $f$. 

Using the fact that the degree of an odd map $S^n\to S^n$ is odd, see~\cite{matousek2003using}, we obtain that the generator $[S^n\times \{*\}]$ of $H_n(S^n\times S^n)$ is mapped by $\widetilde F$ to an odd multiple of the fundamental class $[S^n]$. Since the homology class of the diagonal in $S^n\times S^n$ equals the sum $[S^n\times \{*\}]+[\{*\}\times S^n]$ and $\widetilde F$ maps the diagonal to $\deg f[S^n]\equiv 0[S^n]\mod 2$, the other generator $[\{*\}\times S^n]$ is mapped to an odd multiple of the fundamental class $[S^n]$. Let us now perturb $\widetilde F$ slightly into a smooth map, since we are going to apply Sard's theorem.

Consider the map $S^n\times S^n\times [-1,1]\to S^n\times [-1,1]$ defined by $(x,y,t) \mapsto (\widetilde F(x,y), t)$.
We would like to turn it into a $S^{2n+1}\to S^{n+1}$ map by identifying some points in both domain and codomain. On the top $S^n\times S^n\times \{1\}$ we pick a point $*\in S^n$ and identify $x\times y$ with $x\times *$ for all $x\times y\in S^n\times S^n$. Likewise, on the bottom $S^n\times S^n\times \{-1\}$ we identify $x\times y$ with $*\times y$ for all $x\times y\in S^n\times S^n$. This way the domain becomes the join $S^n * S^n = S^{2n+1}$. In the codomain we identify all the points in $S^n\times \{1\}$ with $*\times \{1\}$ and all the points in $S^n\times \{-1\}$ with $*\times \{-1\}$ so that the codomain becomes $S^{n+1}$. Thus we construct a continuous map $\hat F : S^{2n+1}\to S^{n+1}$. The splitting $S^{n+1} = D_+\cup D_-$ into the northern and southern hemispheres corresponds to the splitting of the sphere $S^{2n+1}$, $S^{2n+1} = T_+\cup T_- = \hat F^{-1}(D_+)\cup \hat F^{-1}(D_-)$, into solid tori.

Let $Z\subseteq S^n\times S^n$ be the $\widetilde F$-preimage of a regular value $y\in S^n$. From the above observation in the $n$-dimensional homology of $S^n\times S^n$ the homology class $Z$ has odd intersection with both generators of $H_n(S^n\times S^n)$. This implies that the $\hat F$-preimage of its regular point $(y,t) \in D_+$ is $Z\times \{t\}$, $1>t>0$, and is homologous to an odd multiple of the $n$-dimensional homology generator of the torus $T_+$. The $\hat F$-preimage of its regular point $(y,t) \in D_-$ is $Z\times \{t\}$, $-1<t<0$, and is homologous to an odd multiple of the $n$-dimensional homology generator of the torus $T_-$. Since the linking number in $S^{2n+1}$ of the $n$-dimensional homology generators of $T_+$ and $T_-$ is $\pm 1$ (the sign depending on the choice of orientation), the linking number in $S^{2n+1}$ of two preimages is odd. 


It remains to apply the Adams theorem \cite{ad1960} asserting that this linking number, the Hopf invariant, can be odd only when $n=1,3,7$.
\end{proof}

\begin{remark}
For $n=1$, a simple exercise on the intermediate value theorem shows that for a continuous periodic function $f :\mathbb R\to \mathbb R$ and any parameter $\delta$ the equation $f(x+\delta)=f(x)$ has a solution, which is a positive result for zero-degree maps $S^1\to S^1$. For other degrees, the maps $z\mapsto z^d$ of the unit circle in the complex plane only glue points at distances $2\pi k / d$ for $k\in\mathbb Z$, which is a negative result. 

In the case $n=3$ one may try explicit formulas with unit quaternions. But the obvious candidate $q\mapsto q^2$ sends the whole imaginary subsphere $S^2\subseteq S^3$ to the single point $-1$, thus giving no counterexample. So the cases $n=3,7$ in Theorem~\ref{theorem:hopf-degree-zero} seem to remain open.
\end{remark}

\section{New Hopf type results}
\label{section:hopf-new}

Now we discuss some new results. We modify the proof of Theorem~\ref{theorem:bu-paths} to obtain the following generalization of the Hopf theorem.

\begin{theorem}
\label{theorem:hopf-paths}
Suppose $X$ is a closed manifold of dimension $n$, $Y$ is an open manifold of dimension $n$, and $f : X\to Y$ is a continuous map. Assume that $X$ has a metric with injectivity  radius $\rho$ and $0<\delta \le \rho$. Then there exist a pair of points $x,y\in X$ at distance $\delta$ such that $f(x) = f(y)$.
\end{theorem}

\begin{remark}
Compared to the Hopf theorem, in this theorem we assume additionally that $\delta$ is at most the injectivity radius, but we allow arbitrary open manifold in place of $\mathbb R^n$ as the codomain.
\end{remark}

\begin{proof}
We mostly follow the proof of Theorem~\ref{theorem:bu-paths}. Assume that $f$ is generic and consider the preimages of a regular value $y\in Y$. Since $Y$ is open, the degree of $f$ is even and $f^{-1}(y)$ consists of an even number of points.

Assuming that no two points in $f^{-1}(y)$ are at distance $\delta$, make a graph $G_y$ on vertices $f^{-1}(y)$ and edges corresponding to pairs $x',x''$ at distance less than $\delta$. By the assumption on the injectivity radius this graph can be drawn by shortest paths on $X$ and depends continuously on $x'$ while $x''$ does not cross special values of $f$.

As in the proof of Theorem~\ref{theorem:bu-paths}, for generic smooth $f : X\to Y$, let $S_1\subset Y$ be the set of special values of $f$, which has codimension at least $1$, and let the set $S_2\subset S_1$ correspond to the singularities of $f$ more complicated than folds, of codimension at least $2$. The space $Y\setminus S_2$ remains connected from dimension considerations whenever $Y$ was connected. Now we see that when the point $y$ travels in $Y\setminus S_2$, the graph $G_y$ may only change when $y$ crosses a fold singularity and some two vertices of the graph are ``created'' or ``annihilated'', having the same sets of neighbors in the remaining vertices of the graph $G_y$.

Now we want to repeat the part of the proof of Theorem~\ref{theorem:bu-paths} using the homotopy in the space of cycles along the edges of $G_y$:
$$
h_t(y) = \sum_{(x_1, x_2)\in E(G_y)} s(x_1, x_2)(t/2).
$$
Like in Remark~\ref{fund-filling}, this homotopy may be interpreted as an $(n+1)$-dimensional chain in $X$. But unlike the proof of Theorem~\ref{theorem:bu-paths}, the mod $2$ boundary of this chain may \emph{not} be the fundamental class $[X]$, but is the set of those points $x\in X\setminus f^{-1}(S_1)$ that come with odd degree in their corresponding graph $G_{f(x)}$. Fortunately, we will show that actually all the points of $X\setminus f^{-1}(S_1)$ have odd degrees in their $G_{f(x)}$ and the proof can be finished similar to the proof of Theorem~\ref{theorem:bu-paths}.

Without loss of generality assume that $X$ is connected and move a point $x$ in $X\setminus f^{-1}(S_2)$, which is also connected from dimension considerations. During such a move there may be two possible modifications of the graph $G_{f(x)}$:

1) a pair of vertices $(x', x'')$ disjoint from $x$ is added or removed from $G_{f(x)}$. Since the points $x'$ and $x''$ have the same sets of neighbors $N(x')\setminus x''=N(x'')\setminus x'$ then the degree of $x$ is changed by $\mp 2$ on this event;

2) the vertex $x$ collides with another vertex $x'$ in $G_{f(x)}$ and they exchange places. Because their sets of neighbors are the same, $N(x)\setminus x' = N(x')\setminus x$, then the degree of $x$ does not change on this event.

Therefore for any $x\in X\setminus f^{-1}(S_1)$ the degree of $x$ in $G_{f(x)}$ is the same mod $2$. Now remember that $Y$ is open and $X$ is closed, then for some $y\in Y$ the graph $G_y$ must be empty and while moving to a nonempty graph it will first generate a pair of points connected by an edge. Hence for some point $x\in X\setminus f^{-1}(S_1)$ its degree in $G_{f(x)}$ must be odd and therefore it must be odd for every $x\in X\setminus f^{-1}(S_1)$. So the image of $h_t(y)$ is a chain in $C_{n+1}(X; \mathbb F_2)$ (see Remark~\ref{fund-filling}) with boundary $[X]$ mod $2$, which is a contradiction, because $X$ is closed.
\end{proof}




Another approach to Hopf type results is possible, following~\cite{vol1992}. Informally, we may increase the dimension of $Y$, drop the compactness assumption on $X$, but require an assumption on its Stiefel--Whitney classes (compare with~\cite[Theorem~1.2]{haibao1989}):

\begin{theorem}
\label{hopf-sw}
Let $f : X\to Y$ be a continuous map between manifolds that induce a zero map on cohomology modulo $2$ in positive dimensions. Suppose $\bar w_k(TX)\neq 0$ (the dual Stiefel--Whitney class), $\dim X + k - 1 \ge \dim Y$, $X$ is a complete Riemannian manifold, and $\delta$ is a prescribed real number. Then there exists a pair $x,y\in X$ such that $f(x) = f(y)$ and the points $x$ and $y$ are connected by a geodesic of length $\delta$.
\end{theorem}

\begin{proof}
Consider the space $S_X$ of pairs $(x, v)$, where $x$ is an arbitrary point in $X$ and $v$ is a unit tangent vector at $x$. This space has a natural $\mathbb Z_2$-action $(x, v)\mapsto (x, -v)$.

For $\mathbb Z_2$-spaces the following invariant is well-known. The natural $\mathbb Z_2$-equivariant map to the one-point space $\pi_{S_X}\to \pt$ induces the map of the equivariant cohomology
$$
\pi_{S_X}^* : H_{\mathbb Z_2}^*(\pt; \mathbb F_2) \to H_{\mathbb Z_2}^*(S_X; \mathbb F_2).
$$
The former algebra $H_{\mathbb Z_2}^*(\pt; \mathbb F_2)=H^*(B\mathbb Z_2; \mathbb F_2) = \mathbb F_2[t]$ is a polynomial ring with one-dimensional generator $t$. The maximal power of $t$ that is mapped nontrivially to the equivariant cohomology of $S_X$ is called the \emph{homological index} of $S_X$ and denoted $\hind S_X$. In~\cite{cf1960} the following is proved: take the maximal $k$ so that the dual Stiefel--Whitney class $\bar w(TX)$ is nonzero, then
$$
\hind S_X = \dim X + k - 1,
$$
under the assumption of this theorem $\hind S_X \ge \dim Y$.

Now consider the map $h : S_X \to X$ defined as follows: start a geodesic from $x$ with tangent $v$ and consider its point $h(x, v)$ at distance $\delta/2$ from $x$. Now the composition $f\circ h$ maps $S_X$ to $Y$ and induces a zero map on the mod $2$ cohomology of positive dimension. By the main result from~\cite{vol1992} we see that some two pairs $(x, v)$ and $(x, -v)$ should be mapped to the same point, which gives the required pair connected by a geodesic of length $\delta$.

\end{proof}

\section{A Borsuk--Ulam--Hopf-type theorem for multivalued maps}

The contents of this section are motivated by a personal discussion with Misha~Gromov. Theorems \ref{theorem:bu-paths} and \ref{theorem:hopf-paths} may be generalized as follows.

\begin{theorem}
\label{hopf-multivalued}
Suppose $X$ is a closed manifold of dimension $n$, $\widetilde X$ is another closed manifold of the same dimension, and $Y$ is an open manifold of the same dimension. Assume that $X$ has a metric with injectivity radius $\rho$ and $0<\delta \le \rho$. Let $g : \widetilde X \to X$ be a map of odd degree and $f : X\to Y$ be a continuous map. Then there exists a pair of points $x,y\in \widetilde X$ such that the distance between $g(x)$ and $g(y)$ is $\delta$ and $f(x) = f(y)$.
\end{theorem}

\begin{proof}
The proof is basically the same as in Theorem \ref{theorem:hopf-paths}. For a generic $y\in Y$, the preimage $f^{-1}(y)$ consists of an even number of points. This map $y\mapsto f^{-1}(y)$ may be extended to non-generic $y$ and considered as a map $f^c : Y\to \cyc_0(\widetilde X)$. As in the previous results, the image of the fundamental class $[Y]$ is nonzero in the homology of the space of cycles. 

Now we consider the map $g_* : \cyc_0(\widetilde X) \to \cyc_0(X)$ that pushes forward the cycles with $g$. The composition $g_*\circ f^c$ is a map $Y \to \cyc_0(X)$ that also has homologically nontrivial $g_*\circ f^c ([Y])$ when the degree of $g$ is odd.

Assuming that for any $y\in Y$ no two points of $g_*\circ f^c(y)$ are at distance $\delta$, we contract the cycle $g_*\circ f^c ([Y])$ in the space of cycles and obtain a contradiction as in the proof of Theorem \ref{theorem:hopf-paths}.
\end{proof}

Going to the limit $\delta\to \pi-0$, we obtain an extension of the Borsuk--Ulam theorem:

\begin{corollary}
\label{bu-multivalued}
Suppose $\widetilde X$ is a closed manifold of dimension $n$ and $Y$ is an open manifold of the same dimension. Let $g : \widetilde X \to \mathbb S^n$ be a map of odd degree and $f : X\to Y$ be a continuous map. Then there exists a pair of points $x,y\in \widetilde X$ such that $g(x)=-g(y)$ and $f(x) = f(y)$.
\end{corollary}

\begin{proof}[A more elementary proof of Corollary~\ref{bu-multivalued} when $Y=\mathbb R^n$]
Consider the map $g\times g : \widetilde X\times \widetilde X \to S^n\times S^n$, this map's degree is the square of the degree of $g$ and is therefore odd. Let 
\[
A = \{(x, -x)\in S^n\times S^n\ |\ x\in S^n\}
\]
be the anti-diagonal. Note that the map $g\times g$ is equivariant with respect to the permutation of factors in the products. Note also that this permutation of factors acts freely on a neighborhood of $A$ and on a neighborhood of $(g\times g)^{-1}(A)$. Then the equivariant version of Thom's transversality theorem applies and we may approximate $g\times g$ by an equivariant $G$ so that $G$ is transversal to $A$.

Now put $M = G^{-1}(A)$, this set is invariant with respect to the permutation of factors in $\widetilde X\times \widetilde X$ and is an $n$-dimensional closed manifold from the transversality. Hence this is a manifold with a free involution. For pairs $(x,y)\in M$, we have
\[
|g(x) + g(y)| < |G_1(x,y) + G_2(x,y)| + 2\varepsilon = 2\varepsilon,
\]
where $G_i$ are the components of $G$ and $G$ approximates $g\times g$ with uniform precision $\varepsilon>0$. Now we restrict the map $(x,y)\mapsto f(x) - f(y)$ to $M$ thus obtaining the map $F : M\to \mathbb R^n$. This map is an equivariant map $F: M\to \mathbb R^n$ (with respect to the antipodal involution of $\mathbb R^n$) and we need to show that $F$ maps some pair $(x,y)\in M$ to zero.

The solution of $F(x,y)$ is guaranteed when $M$, as a smooth closed manifold with involution, has the Borsuk--Ulam property, see \cite{musin2012} for a thorough investigation of this situation. We will show this property by considering another equivariant map $\Phi : M\to \mathbb R^n$ in place of $F$, which is transversal to zero and has an odd number of orbits (of the involution) going to zero. Then the zero set of $F$ must also be non-empty, since a generic smooth equivariant homotopy of $\Phi$ to $F$ establishes an equivariant bordism between the solution sets $\{F(x)=0\}$ and $\{\Phi(x)=0\}$, keeping the parity of the number of orbits in those finite sets.

In order to build an equivariant $\Phi : M\to \mathbb R^n$ we consider the standard projection $P : S^n\to \mathbb R^n$ that is transversal to zero and has precisely one pair of points (an orbit of the antipodal involution) in the preimage of zero. Projecting in an appropriate direction, we may choose this pair $P^{-1}(0)$ to be a regular value of $G_1$. Then the composition $\Phi=P\circ G_1$ is also equivariant (since $G_1(x,y) = - G_2(x,y) = -G_1(y,x)$) and has an odd number of orbits of the involution in $\Phi^{-1}(0)$. This establishes the Borsuk--Ulam property for $M$ and shows that $F^{-1}(0)$ is also non-empty.

Hence we have found $(x,y)\in \widetilde X\times \widetilde X$ such that $f(x) = f(y)$ and $|g(x) + g(y)| < 2\varepsilon$. Then we take the limit $\varepsilon \to +0$ and use the compactness of $\widetilde X\times \widetilde X$ to find a precise solution such that
\[
f(x) = f(y)\quad\text{and}\quad g(x) = - g(y).
\]
\end{proof}

These results imply the following necessary and sufficient conditions for nonzero $1$-Lipschitz maps between ellipsoids.

\begin{theorem}
\label{ellipsoids-lipschitz}
Let $0<a_1\le \dots \le a_n$ and $0<b_1\le \dots \le b_n$ be sequences of reals and let
\[
E_a = \left\{ \frac{x_1^2}{a_1^2} + \dots + \frac{x_n^2}{a_n^2} =1 \right\}
\quad\text{and}\quad
E_b = \left\{ \frac{x_1^2}{b_1^2} + \dots + \frac{x_n^2}{b_n^2} =1 \right\}
\]
be (surfaces of) ellipsoids in $\mathbb R^n$. If there exists a $1$-Lipschitz (in the extrinsic metric of $\mathbb R^{n+1}$) odd degree map $L : E_a\to E_b$ then the inequalities $a_k\ge b_k$ hold for every $k$.
\end{theorem}

\begin{proof}
Assume the contrary, that $a_k < b_k$ for some $k$. Let $X$ be the subset of $E_b$ satisfying $x_1=\dots = x_{k-1} = 0$, this is another (surface of an) ellipsoid of dimension $n-k$ with all axes greater or equal to $b_k$. We first perturb $L$ making it smooth and increasing its Lipschitz constant by an arbitrarily small amount. Then using Thom's transversality theorem we further perturb $L$ so that $\widetilde X = L^{-1}(X)$ is a submanifold of $E_a$. This may again spoil the Lipschitz constant by arbitrarily small amount, but the $1$-Lipschitz property may be compensated by a slight inflation of $E_a$ keeping the assumption $a_k < b_k$.

Let $f : \widetilde X\to \mathbb R^{n-k}$ be the projection setting the first $k$ coordinates to zero. And let $f : \widetilde X \to X$ be the restriction of $L$, the definition of the mapping degree and the transversality assumption guarantee that $f$ has nonzero degree. Applying Theorem \ref{bu-multivalued} to this situation we find two points $x,y\in E_a$ such that their images $f(x)=L(x), f(y)=L(y)\in X$ are opposite and their last $n-k$ coordinates are equal. The former property ensures $|L(x) - L(y)| \ge 2b_k$, while the latter implies $|x-y|\le 2a_k$. Together with $a_k<b_k$ this contradicts the $1$-Lipshitz assumption. 
\end{proof}

\begin{remark}
Theorem \ref{ellipsoids-lipschitz} fails for the intrinsic metrics of $E_a$ and $E_b$ as stated. Hence there remains an open question what kind of claim is true for the two surfaces of an ellipsoid considered with their intrinsic metrics.
\end{remark}

\begin{remark}
If we replace the assumption on odd degree of $L$ by the assumption that $L$ is a homeomorphism in Theorem \ref{ellipsoids-lipschitz} then the proof passes with the classical Borsuk--Ulam theorem in place of Corollary \ref{bu-multivalued}, since in this case $\widetilde X \cong X$.
\end{remark}

\begin{remark}
Our argument in the proof of Theorem \ref{ellipsoids-lipschitz} does not pass for maps of non-zero even degree. At the moment we do not know if this theorem holds in this case.
\end{remark}

\bibliography{../Bib/karasev}
\bibliographystyle{abbrv}
\end{document}